\documentclass[11pt,a4paper]{article}

\usepackage[T1]{fontenc}
\usepackage{lmodern}
\usepackage{microtype}
\usepackage{amsmath,amssymb,amsthm,color}
\usepackage[margin=1in]{geometry}
\usepackage[hidelinks]{hyperref}
\usepackage{authblk}
\usepackage{cite}
\newtheorem{theorem}{Theorem}[section]
\newtheorem{lemma}[theorem]{Lemma}

\newtheorem{conjecture}[theorem]{Conjecture}
\newtheorem{cor}[theorem]{Corollary}
\newtheorem{fact}{Fact}

\theoremstyle{definition}
\newtheorem{remark}{Remark}

\newenvironment{wst}
{\setlength{\leftmargini}{1.5\parindent}
	\begin{itemize}
		\setlength{\itemsep}{-1.1mm}}
	{\end{itemize}}

\begin{document}
	
	\title{A degree-sum condition for paths through prescribed vertices\thanks{Leilei Zhang's work was supported by  the National Natural Science Foundation of China (Grant Nos. 12601664, 12571365), JSPS KAKENHI (Grant No.~25KF0036), the Hubei Provincial Natural Science Foundation of China (Grant No.~2025AFB309), the China Postdoctoral Science Foundation (Grant No.~2025M773113), the Fundamental Research Funds for the Central Universities, Central China Normal University (Grant No.~CCNU24XJ026).}}
	
	\author[1]{Kun Cheng\footnote{Email: \texttt{chengkunmath@163.com}}}
	\author[2]{Leilei Zhang\footnote{Email: \texttt{mathdzhang@163.com}}}
	
	\date{}
	
	\affil[1]{\small Department of Mathematics, East China Normal University,  Shanghai, 200241, China }
	\affil[2]{\footnotesize School of Mathematics and Statistics, Key Laboratory of Nonlinear Analysis \& Applications (Ministry of Education), Central China Normal University, Wuhan 430079, China}
	
	\maketitle
	
	\begin{abstract}
		Let $G$ be a graph of order $n$ and let $S\subseteq V(G)$. Denote by $\kappa(S)$ the minimum size of a vertex set separating two vertices of $S$, and by $\sigma_p(S)$ the minimum degree sum over all independent sets of $p$ vertices in $S$. We prove that if $\kappa(S)\ge k\ge2$ and $\sigma_p(S)\ge n+(p-1)^2$ for every integer $p\ge k$, then any two distinct vertices of $S$ are joined by a path containing all vertices of $S$. This gives a path version of Ota's theorem on cycles containing a prescribed set of vertices. As an application, we prove a conjecture of Faudree, Faudree, Gould, Jacobson and Magnant: for $k\ge3$, every $k$-connected graph $G$ of order $n$ satisfying $\delta(G)>(n+k^2-2k)/k$ and $\delta(G)\ge\alpha(G)+k-2$ is hamiltonian-connected.
	\end{abstract}
	
	\noindent\textbf{Keywords:} Hamiltonian-connected;  Degree sum; Connectivity; Independence number.
	
	\noindent\textbf{2020 Mathematics Subject Classification:} 05C45,
	05C38, 05C40.
	
	\section{Introduction}
	We consider finite simple graphs and use standard terminology and notation from \cite{BM} and \cite{W}.    
	For a graph $G=(V(G),E(G))$, let $\delta(G)$, $\kappa(G)$, and $\alpha(G)$ denote its minimum degree, connectivity, and independence number, respectively.  For $x\in V(G)$, write $N_G(x)$ for the neighborhood of $x$ and $d_G(x)=|N_G(x)|$ for its degree. 
	For two distinct vertices $x$ and $y$ in $G$, an {\em $(x, y)$-path}  is a path whose endpoints are $x$ and $y$.  
	A Hamilton cycle (resp., path) in $G$ is a cycle (resp., path) containing every vertex of $G$.  
	A graph is {\em hamiltonian} if it contains a Hamilton cycle, and it is {\em hamiltonian-connected} if every two distinct vertices are joined by a Hamilton path.

	The existence of Hamilton cycles is a classical topic in graph theory. 
	Dirac \cite{D} proved that every graph $G$ of order $n\ge3$ with $\delta(G)\ge n/2$ is hamiltonian.   
	Ore \cite{OO} replaced Dirac's condition by $d_G(x)+d_G(y)\ge n$ for every pair of nonadjacent vertices $x$ and $y$.  
	A well-known result of Chv\'{a}tal and Erd\H{o}s~\cite{CE} from 1972 states that a graph whose connectivity is at least its independence number is hamiltonian, and that a graph whose connectivity is at least one more than its independence number is hamiltonian-connected.

	Thus degree, connectivity, and independence number provide different ways to ensure a Hamilton cycle or a Hamilton path.  See
	\cite{JO} for a survey of path and cycle results based on connectivity and independence number.     
    For a $k$-connected graph, the Chv\'{a}tal–Erd\H{o}s theorem~\cite{CE} ensures Hamiltonicity when $\alpha(G)\le k$. Under an additional minimum-degree condition, Faudree, Faudree, Gould, Jacobson and Magnant~\cite[Proposition 1]{FFGJM} replaced this bound by $\alpha(G)\le\delta(G)-k+2.$ They proved that:
	If $G$ is a $k$-connected graph of order $n$ with 
	$\delta(G)>(n+k^2-k-1)/(k+1)$
	and 
	$\delta(G)\ge\alpha(G)+k-2,$  where $k\ge2$, 
	then $G$ is hamiltonian. 
	To study the hamiltonian-connectedness, they proposed the following conjecture.
	\begin{conjecture}[Faudree--Faudree--Gould--Jacobson--Magnant \cite{FFGJM}]\label{conj:ffgjm}
		Let $G$ be a $k$-connected graph of order $n$, where $k\ge3$.  If
		$\delta(G)>(n+k^2-2k)/k$ 
		and $\delta(G)\ge\alpha(G)+k-2,$ 
		then $G$ is hamiltonian-connected.
	\end{conjecture}

	As a first step, they showed that Conjecture~\ref{conj:ffgjm} holds in the following two cases.
	
	\begin{theorem}[Faudree--Faudree--Gould--Jacobson--Magnant \cite{FFGJM}] 
		Let $k\ge3,$ and let $G$ be a graph of sufficiently large order $n$
		satisfying
		$
		\delta(G)>(n+k^2-2k)/k
		$ and 
		$\delta(G)\ge\alpha(G)+k-2.$
		Then $G$ is hamiltonian-connected if one of the following holds:
		\begin{wst}
			\item[{\rm (i)}] $\kappa(G)\ge4k^2+1$;
			\item[{\rm (ii)}] $k\in\{3,4\}$ and $\kappa(G)\ge k$.
		\end{wst}
	\end{theorem}
	
	Thus the cases $k=3,4$ for Conjecture~\ref{conj:ffgjm} were known for sufficiently large graphs, while the
	general result required $\kappa(G)\ge 4k^2+1.$   
	The gap $k\le \kappa(G)<4k^2+1$ and the restriction on the order were left open.
	
	There is also a prescribed-vertex version of the hamiltonian problem.  A set $S\subseteq V(G)$ is called \emph{cyclable} if there exists one cycle of $G$ containing every vertex of $S$.  Dirac's  theorem~\cite{GAD} states that every $k$ specified vertices of a $k$-connected graph lie on a common cycle. It led to many later results on paths and cycles through specified vertices; see, for example, \cite{EGL,FLMW,G,H}.  
	
	Let $G$ be a graph and let $S\subseteq V(G)$ with $|S|\ge2$.  A set $U\subseteq V(G)\setminus\{x,y\}$ \emph{separates} distinct vertices $x,y$ if $x$ and $y$ lie in different components of $G-U.$  If $G[S]$ is not complete, define
	$$
	\kappa(S)=\min\bigl\{|U|: U\text{ separates two vertices of }S\bigr\}.
	$$
	If $G[S]$ is complete, define $\kappa(S)=|S|-1$.  In particular, $\kappa(V(G))=\kappa(G)$, where $\kappa(K_n)=n-1$.  Write
	$\alpha(S)$ for the maximum cardinality of an independent subset of $S$. For $1\le p\le\alpha(S)$, define
	$$
	\sigma_p(S)= \min\left\{  \sum_{x\in X}d_G(x): X\subseteq S,~X\text{ is independent},~|X|=p  \right\}.
	$$
	Write $\sigma_p(G)=\sigma_p(V(G))$ for short.  
	In 1995, Ota obtained the following degree-sum condition for cyclability.
	
	\begin{theorem}[Ota \cite{KO}]\label{thmOta}
		Let $G$ be a graph of order $n$ and let $S\subseteq V(G)$ with $\kappa(S)\geq k \ge2$.  If
		$\sigma_p(S)\ge n+(p-1)(p-2)$ 
		for every integer $p\ge k+1$,
		then $S$ is cyclable.
	\end{theorem}

	Theorem~\ref{thmOta} has hamiltonian consequences and led to later degree-sum results; see \cite{BLLTV,OY}.  
    Our aim is to prove a path version in which two endpoints are prescribed. 
	Long paths between prescribed endpoints have also been studied under degree and Chv\'atal--Erd\H{o}s type conditions.  
    Hu, Tian, Wei, Egawa and Hirohata \cite{HTWEH} proved a Fan-type lower bound for a path joining any prescribed pair of vertices in a $3$-connected graph.  Chen, Hu and Wu \cite{CHW} proved a related Chv\'atal--Erd\H{o}s type result. Those results use hypotheses on the graph and seek long paths between two arbitrary endpoints.  Our main result uses degree sums  for independent subsets of $S$ and requires a path passing through all vertices of $S$.

	\begin{theorem}\label{thm:main}
		Let $G$ be a graph of order $n$ and let $S\subseteq V(G)$ with $|S|\ge2$.
		If $\kappa(S)\ge k\ge2$ and
		\begin{align}\label{eqCnZ1}
			\text{
				$\sigma_p(S)\ge n+(p-1)^2$ for every integer $k\le  p\le \alpha(G),$ 
			}
		\end{align}
		then for every two vertices $u,v\in S$, there is a $(u,v)$-path containing all vertices of $S$.
	\end{theorem}
	
	A non-complete graph $G$ is said to be {\em $t$-tough} if for every vertex cut $S$ of $G$, the ratio of $|S|$ to the number of components of $G-S$ is at least $t$. The {\em toughness} $\tau(G)$ of  $G$ is the maximum value of $t$ such that $G$ is $t$-tough. 
    (taking $\tau(K_n)=\infty$ for
	all $n\ge 1$). The concept of toughness of a graph was introduced by
	Chv\'{a}tal~\cite{Chvatal}. 
    Clearly, every hamiltonian-connected has toughness greater than one. 
    
    Substituting $S=V(G)$ in Theorem~\ref{thm:main} gives the following. 
	\begin{cor}\label{cor31}
		Let $G$ be a $k$-connected graph, where $k\ge 2.$ If 
		\begin{align}\label{eqcz3}
			\text{
				$\sigma_p(G)\ge n+(p-1)^2$ for every integer $k\le p\le \alpha(G),$ 
			}
		\end{align}
		then $G$ is hamiltonian-connected.  
	\end{cor}
	
	\begin{remark}\label{r1}
		The condition in \eqref{eqcz3}, and hence \eqref{eqCnZ1}, can not be weakened.  
		In fact, if we replace \eqref{eqcz3} with 
		\begin{align}\label{eq2}
			\text{
				$\sigma_p(G)\ge n+p^2-2p$~~for every integer $k\le p\le \alpha(G),$ 
			}
		\end{align}
		then let $H=K_k\vee kK_m.$ 
		Clearly, 
		$\alpha(H)=\kappa(H)=k.$ Thus there exists only one integer $p$ satisfying $k\le p\le \alpha(G),$ and 
		$\sigma_k(H)=k\delta(H)=k(k+m-1)=|H|+k^2-2k.$ 
		Hence $H$ satisfies \eqref{eq2}. 
		However, $H$ is not hamiltonian-connected, because its toughness is at most 1.
	\end{remark}
	
	As an application of Theorem \ref{thm:main}, we prove Conjecture \ref{conj:ffgjm}. We organize the paper as follows. In Section~\ref{s2cz}, we present  the proof of Theorems~\ref{thm:main}.  In Section~\ref{s3cz} we prove Conjecture~\ref{conj:ffgjm}.
	
	\section{Proof of Theorem~\ref{thm:main}}\label{s2cz}
	Let $G$ be a graph.
	For a vertex subset $Y\subseteq V(G)$ and a vertex $x\in V(G),$ let $N_Y(x)=N_G(x)\cap Y$ and let $d_Y(x)=|N_Y(x)|.$
	Let $P=v_1v_2\cdots v_\ell$ be a path in a graph $G.$ 
	We set 
	\begin{align*}
		\text{
			$v_i^+:=v_{i+1}$~~and~~ 
			$v_j^-:=v_{j-1}$,
	}\end{align*}
	respectively.  
	If $a<b,$
	we define two subpaths 
	\begin{align*}
		\overrightarrow{P}[v_a,v_b]:=v_a v_{a+1}\cdots  v_b \text{~~and~~} 
		\overleftarrow{P}[v_b,v_a]:=v_b v_{b-1}\cdots  v_a.
	\end{align*}   
	The notation $\overrightarrow{P}(v_a,v_b]$, $\overrightarrow{P}[v_a,v_b)$, and $\overrightarrow{P}(v_a,v_b)$ means that the endvertex $a,$ the endvertex $b,$ or both endvertices $\{a,b\},$ respectively, are deleted from $\overrightarrow{P}[v_a,v_b].$  
	We also allow $\overrightarrow{P}[v_a,v_a]$ to be a trivial path; every open or half-open version with equal endpoints is empty. 
	Let 
	$$v_a\overrightarrow{P}v_b=\overrightarrow{P}[v_a,v_b]~~\text{and}~~
	v_b\overleftarrow{P}v_a=\overleftarrow{P}[v_b,v_a].$$

	A family of $k$ internally disjoint $(x,Y)$-paths whose terminal vertices are distinct is referred to as a {\em $k$-fan} from $x$ to $Y$. 
	We denote by $f(x,Y)$ the maximum integer $k^*$ such that there exists a $k^*$-fan from $x$ to $Y$.
	\begin{lemma}\label{lem:fan}
		Let $G$ be a graph and let $S\subseteq V(G)$ with $\kappa(S)\ge k$. Then 
		for every vertex $z\in S,$ 
		$f(z, S\setminus\{z\})\ge \min\{k,|S|-1\}.$  
	\end{lemma} 
	\begin{proof} 
		If $S$ is a clique of $G,$ 
		then for any $z\in S,$ 
		$z$ is adjacent to every vertex of $S\setminus\{z\}.$ 
		Thus there exists a $(|S|-1)$-fan from $z$ to $S\setminus \{z\},$ 
		as desired. 
		
		Next we assume that $G[S]$ is not a complete graph. 
		Let $G'$ be a graph obtain from $G$ by  adding a new vertex $w$ and joining $w$ to every vertex of $S\setminus\{z\}.$ 
		Suppose to the contrary that there exists no $\min\{k,|S|-1\}$-fan in $G$ 
		from $z$ to $S\setminus\{z\},$ for some vertex $z\in S.$ 
		Then Menger's theorem~\cite{M} implies that there exists a separating set $U\subseteq V(G')\setminus\{z,w\}$ with $|U|\le \min\{k,|S|-1\}-1$ which separates $z$ and $w.$ 
		Hence, in $G-U,$ there exists a vertex in $S\setminus (U\cup \{z\})$ that belongs to a component different from the one containing $z.$ Therefore, 
		$\kappa(S)\le |U|\le k-1,$ a contradiction. 
	\end{proof}

	Next we prove Theorem~\ref{thm:main}. 
	\begin{proof}[\bf Proof of Theorem~\ref{thm:main}]
		Let $G$ be a  graph of order $n,$  
		and let $S\subseteq V(G)$ with $|S|\ge 2.$  
		Assume that $\kappa(S)\ge k\ge 2$  
		and 
		$\sigma_p (S)\ge n+p^2-2p+1$ for any integer $p$ with $k\le p\le \alpha(G).$ 
		
		We use induction on $|S|.$ 
		If $|S|=2,$ 
		then by $\kappa(S)\ge k\ge 2,$  
		we are done. 
		Next we assume that $|S|\ge 3$ and the result holds for every vertex subset $T$ with $2\le|T|\le |S|-1.$ 	
		Now we are to show that for any two vertices $u',v'\in S,$ 
		there exists an $(u',v')$-path in $G$ containing all vertices of  $S.$ 
		
		Suppose to the contrary that there are two vertices $u,v\in S$
		for which no $(u,v)$-path in $G$ contains all vertices of $S$.
		We choose the vertex $w\in S\setminus\{u,v\}$ such that $f(w,S\setminus\{w\})$ is as large as possible. 
		Note that 
		\begin{align*}
			\text{
				$\kappa\left(S\setminus\{w\}\right)\ge \kappa(S)\ge k$~~and~~$\sigma_p(S\setminus\{w\})\ge \sigma_p(S)\ge n+p^2-2p+1$ for any $p\ge k.$ 
			}
		\end{align*}
		Then by the induction hypothesis, 
		there exists a $(u,v)$-path  $P$  in $G$ such that 
		$S\setminus \{w\}\subseteq V(P).$ 
		Clearly, $w\not\in V(P).$ 
		Based on $\kappa(S)\ge k$  and Lemma~\ref{lem:fan},  we have that 
		there exist $q :=f(w,S\setminus\{w\})\ge \min\{k,|S|-1\}$ paths $P_1,\ldots, P_q$ from $w$ to $V(P)$ such that 
		$V(P_i)\cap V(P_j)=\{w\}$ for $i\neq j.$  	
		Let $a_i$ be the vertex of $V(P_i)\cap V(P),$ and assume that $a_1,\ldots, a_q$ appear in this order along $\overrightarrow{P}$. 	
		We choose $P_1,\ldots, P_q$ such that
		$\sum_{i=2}^q \left|  \overrightarrow{P}[a_1,a_i] \right|$ is minimum.

		\begin{fact}\label{f1} 
			$V\left( \overrightarrow{P}(a_{i-1},a_i)\right)\cap S\neq \emptyset$ 
			for any $2\le i\le q.$ 
		\end{fact} 
		\begin{proof}[Proof of Fact~\ref{f1}]
			Suppose to the contrary that there exists an integer $i$  such that 
			$V\left( \overrightarrow{P}(a_{i-1},a_i) \right)\cap S= \emptyset$  and $2\le i \le q.$ 
			Then
			$$
			u \overrightarrow{P}a_{i-1} P_{i-1} w P_i a_i \overrightarrow{P}v
			$$
			is a $(u,v)$-path in $G$ passing through all vertices of $S,$ a contradiction. This proves Fact~\ref{f1}.  
		\end{proof} 
		
		By Fact~\ref{f1}, for each $2\le i\le q,$   let $x_i$ denote the last vertex of 
		$\overrightarrow{P}(a_{i-1},a_i)$ that belongs to $S.$ 
		We define $x_1$ as follows:
		if $u\neq a_1,$ 
		let $x_1$ be the last vertex of $V\left( \overrightarrow{P}[u,a_1) \right)\cap S;$  
		if $u= a_1,$ let $x_1=u.$ 
		\begin{fact}\label{f2}
			\begin{wst}
				\item[{\rm (i)}] $|S|\ge k+1,$ and hence $q\ge  k\ge 2;$
				\item[{\rm (ii)}]
				$N(x_j)\cap V\left( \overrightarrow{P}[x_i,a_i) \right)=\emptyset$ for any two distinct integers $i,j$ with $2\le j\le q$ and $1\le i\le q.$  
			\end{wst}
		\end{fact}
		\begin{proof}[Proof of Fact~\ref{f2}]
			(i) Suppose to the contrary that $|S|\le k.$ 
			Since $\{x_2,\ldots, x_q\}\cup \{u,v,w\}\subseteq S,$ 
			we have  $|S|\ge q+2,$ 
			which contradicts the fact that $q\ge \min\{k,|S|-1\}=|S|-1.$ 
			This implies that $|S|\ge k+1,$ 
			and thus $q\ge \min\{k,|S|-1\}= k.$

			(ii)
			Suppose to the contrary that there exist two distinct integers $i,j$ such that $2\le j\le q,$ $1\le i\le q$ and $N(x_j)\cap V\left( \overrightarrow{P}[x_i,a_i) \right)\neq \emptyset.$ Let $z\in N(x_j)\cap V\left( \overrightarrow{P}[x_i,a_i) \right).$ 
			If $j<i,$ then 
			$$u \overrightarrow{P}x_j z \overleftarrow{P} a_j P_j w P_i a_i \overrightarrow{P} v$$ 
			is a $(u,v)$-path passing through all vertices of $S,$ a contradiction. 
			Thus we assume that $i<j.$ 
			Then 
			$$u \overrightarrow{P} z x_j \overleftarrow{P} a_i P_i w P_j a_j \overrightarrow{P} v$$ 
			is a $(u,v)$-path containing all vertices of $S$, yielding a contradiction. 
			This proves Fact~\ref{f2}. 
		\end{proof} 
		
		\begin{fact}\label{f3}
			$\{w\}\cup \{x_2,\ldots, x_q \}$ is an independent set of $G.$ 
		\end{fact}
		\begin{proof}[Proof of Fact~\ref{f3}] 
			By Fact~\ref{f2}, $x_i\not\sim x_j$ for any two distinct  $2\le i,j \le q.$ 
			Suppose that there exists an integer $i\in\{2,\ldots,q\}$ such that  
			$w\sim x_i.$ 
			Then 
			$u\overrightarrow{P} x_i w P_i a_i \overrightarrow{P} v$ 
			is a $(u,v)$-path passing through all vertices of $S,$ a contradiction.  
			Therefore $w\not\sim x_j$ for any $2\le j\le q.$ This proves Fact~\ref{f3}.
		\end{proof} 
		
		\begin{fact}\label{f4}
			\begin{wst}
				\item[{\rm (i)}] 
				$N(x_j)\cap \left( V(P_i)\setminus \{a_i\} \right)=\emptyset$ 
				for any two distinct integers $i,j$ with $1\le i\le q$ and $2\le j\le q;$ 
				\item[{\rm (ii)}] 
				$N(x_i)\cap \left(V(P_i)\setminus \{a_i\} \right)=\emptyset$ for any $2\le i\le q.$ 
			\end{wst}
		\end{fact}
		
		\begin{proof}
			[Proof of Fact~\ref{f4}] 
			(i) 
			Suppose to the contrary that there exist two distinct integers $i,j$ such that $N(x_j)\cap \left( V(P_i)\setminus \{a_i\} \right)\neq \emptyset.$ 
			Let $z\in N(x_j)\cap \left( V(P_i)\setminus \{a_i\} \right).$ 
			If $j<i,$ 
			then 
			$$u\overrightarrow{P}x_j z P_i w P_j a_j \overrightarrow{P} v$$ 
			is a $(u,v)$-path passing through all vertices of $S,$ a contradiction.  Next we assume that 
			$i<j.$ 
			Then 
			$$u \overrightarrow{P} x_j z P_i w P_j a_j \overrightarrow{P} v$$ 
			is a $(u,v)$-path containing all vertices of $S$, yielding a contradiction. 
			Thus (i) holds. 
			
			(ii) 
			Suppose to the contrary that there exists an integer $i\in\{2,\ldots, q\}$  such that 
			$N(x_i)\cap \left(V(P_i)\setminus \{a_i\} \right)\neq \emptyset.$ 
			Let $z\in N(x_i)\cap \left(V(P_i)\setminus \{a_i\} \right).$ 
			Then the path $w P_i z x_i$ contradicts  the choice of 
			$P_1,\ldots,P_q.$ This proves Fact~\ref{f4}. 
		\end{proof} 
		We write $A:=\{a_1,\ldots,a_q\},$ 
		$B:=V(P_1)\setminus \{a_1\}$ 
		and $R:=V(G)\setminus \left( V(P)\cup V(P_1)  \right).$ 
		\begin{fact}\label{f6} 
			\begin{wst}
				\item[{\rm (i)}] 
				$d_B(w)+\sum_{j=2}^q d_B(x_j)\le |B|-1;$ 
				\item[{\rm (ii)}] 
				$d_{V(P)}(w)\le q.$
			\end{wst}
		\end{fact}
		\begin{proof}[Proof of Fact~\ref{f6}]
			(i)
			From Fact~\ref{f4}(i), we have 
			$N(x_j)\cap \left( V(P_1)\setminus \{a_1\}  \right)=\emptyset$ 
			for each $2\le j\le q.$ 
			Therefore, $\sum_{j=2}^q d_B(x_j)=0.$ 
			This, together with the fact that $N_B(w)\subseteq B\setminus\{w\},$ implies $d_B(w)+\sum_{j=2}^q d_B(x_j)\le |B|-1,$ 
			as desired. 
			
			(ii) By the choice of $w$ and $P_1,\ldots, P_q,$ 
			we have  $N(w)\cap V(P)\subseteq \{a_1,\ldots, a_q\}$ and hence $d_{V(P)}(w)\le q.$ 
			This proves Fact~\ref{f6}.  
		\end{proof}

		\begin{fact}\label{f66}
			\begin{wst}
				\item[{\rm (i)}] 
				$N(w)\cap N(x_i)\cap R=\emptyset$ for any integer $i\in \{2,\ldots,q\};$ 
				\item[{\rm (ii)}] 
				$N(x_j)\cap N(x_i)\cap R=\emptyset$ 
				for any two distinct $i,j\in\{2,\ldots, q\};$ 
				\item[{\rm (iii)}] 
				$d_R(w)+\sum_{i=2}^q d_R(x_i)\le |R|.$
			\end{wst}
		\end{fact} 
		
		\begin{proof}
			[Proof of Fact~\ref{f66}] 
			(i) 
			Suppose to the contrary that there exists a vertex $z\in N(x_i)\cap R$ and $z\sim w,$ where $i\in\{2,\ldots, q\}.$ Then by Fact~\ref{f4}(i)-(ii), we have that $z\not\in \bigcup_{i=1} ^q V(P_i).$  
			Hence, 
			$$u\overrightarrow{P} x_i z w P_i a_i \overrightarrow{P} v$$ 
			is a $(u,v)$-path passing through all vertices of $S,$ a contradiction. 
			
			(ii) Suppose to the contrary that there exists a vertex $z\in N(x_i)\cap R$ and $z\sim x_j,$ where $i,j\in\{2,\ldots, q\}$ and $i\neq j.$ 
			Without loss of generality, we assume that  $j<i.$ 
			Then by Fact~\ref{f4}(i)-(ii), we have that
			$$u\overrightarrow{P} x_j z x_i \overleftarrow{P} a_j P_j w P_i a_i  \overrightarrow{P} v$$ 
			is a $(u,v)$-path passing through all vertices of $S,$ a contradiction. 
			
			(iii) From (i) and (ii), we have that 
			\begin{align*}
				d_R(w)+\sum_{i=2}^q d_R(x_i)=\left|(N(w)\cap R)\cup \left( \bigcup_{i=2}^q (N(x_i)\cap R )\right)\right|\le |R|.
			\end{align*} 
			This completes the proof of Fact~\ref{f66}.
		\end{proof}

		A  subpath $P^*=\overrightarrow{P}[z,z')$ 
		of  $P$ is called {\em a good subpath} 
		if it satisfies the following: 
		\begin{wst}
			\item[{\rm (1)}] 
			$z\in S\cup A$ and $z' \in S\setminus A;$
			\item[{\rm (2)}]  
			$V\left( \overrightarrow{P}(z,z') \right) \cap \left( S\setminus A\right)=\emptyset.$ 
		\end{wst}         
		Write $$\mathcal{P}=\{v\}\cup \left\{P^*: P^*\text{ is a good subpath of $P$}\right\}.$$  
		Hence, 
		$$ 
		|\mathcal{P}|\ge q-1 ~~\text{and}~~
		V(P)=\bigcup_{i=1}^q V\left( \overrightarrow{P}[ x_i,a_i) \right) \cup \bigcup_{P^*\in \mathcal{P}} V(P^*).$$ 
		
		\begin{fact}\label{f88} 
			Let $P^*\in \mathcal{P}\setminus\{v\},$ and  
			assume that $P^*=\overrightarrow{P}[z,z')\subseteq \overrightarrow{P}[a_s, a_{s+1}],$ 
			where $1\le s\le q-1.$ Let $g\in V(P^*).$
			The following properties hold:
			\begin{wst}
				\item[{\rm (i)}] 
				If $g\sim x_i$ for some $i\ge s+1,$ 
				then for any $j\ge i+1,$ $N(x_j)\cap V\left( \overrightarrow{P}(g,z') \right)=\emptyset;$
				\item[{\rm (ii)}] If $g\sim x_i$ for some $i\le s,$ 
				then for any $i+1\le j\le s,$ 
				$N(x_j)\cap V\left( \overrightarrow{P}(g,z') \right)=\emptyset;$
				\item[{\rm (iii)}] 
				If $g\sim x_j$ for some $j\ge s+1,$ 
				then for any $i\le s,$ 
				$N(x_i)\cap V\left( \overrightarrow{P}(g,z') \right)=\emptyset.$
			\end{wst}
		\end{fact} 
		\begin{proof}[Proof of Fact~\ref{f88}] 
			(i) 
			Suppose to the contrary that there exists a vertex $z''$ in $N(x_j)\cap V\left( \overrightarrow{P}(g,z') \right)$ for some $j\ge i+1.$ Then $$u\overrightarrow{P} g x_i \overleftarrow{P} z'' x_j \overleftarrow{P} a_i P_i w P_j a_j \overrightarrow{P} v$$ 
			is a $(u,v)$-path passing through all vertices of $S,$ a contradiction.

			(ii) 
			Suppose to the contrary that there exists a vertex $z''$ in $N(x_j)\cap V\left( \overrightarrow{P}(g,z') \right)$ for some $i+1\le j\le s.$ Then $$u\overrightarrow{P} x_i g \overleftarrow{P} a_j P_j w P_i a_i\overrightarrow{P} x_j z'' \overrightarrow{P} v$$ 
			is a $(u,v)$-path passing through all vertices of $S,$ a contradiction. 
			
			(iii)
			Suppose to the contrary that there exists a vertex $z''$ in $N(x_i)\cap V\left( \overrightarrow{P}(g,z') \right)$ for some $i\le s.$ 
			Then $$u\overrightarrow{P} x_i z'' \overrightarrow{P} x_j g \overleftarrow{P} a_i P_i w P_j a_j \overrightarrow{P} v$$ is a $(u,v)$-path passing through all vertices of $S,$ a contradiction.  This proves Fact~\ref{f88}. 
		\end{proof}
		Let $X:=\{x_2,\ldots,x_q\}.$
		For each $P^*\in \mathcal{P},$  let $ |N_X(P^*)|$ be the number of vertices in $X$ that is adjacent to a vertex of $P^*.$  
		\begin{fact}\label{f7}
			$\sum_{i=2}^q d_{V(P^*)}(x_i)\le |N_X(P^*)| +|P^*| -1,$ 
			for each $P^*\in \mathcal{P}.$ 
		\end{fact} 
		\begin{proof}[Proof of Fact~\ref{f7}]
			If $P^*=v,$ 
			then $|P^*|=1$ and $\sum_{i=2}^q d_{V(P^*)} (x_i)=   |N_X(P^*)| +|P^*| -1,$  as desired. 
			
			Next assume that $P^*\neq v.$  
			Let $\ell=|P^*|$ and let  $z_1,\ldots, z_\ell$ be the vertices of $P^*$ in the order on $\overrightarrow{P}$.  
			If $j\ge 2$ and $N_X(v_j)\neq \emptyset$, let $x_{i_j}$ 
			denote the vertex in $N_X(v_j)$ such that $i_j$ is maximum. 
			Define 
			\begin{align*}
				N'_X(v_j)
				&=\begin{cases} N_X(v_j)\setminus\{x_{i_j}\} &{\rm if}\,\,\,j\ge 2\,\,\,{\rm and}\,\,\,N_X(v_j)\neq \emptyset,\\
					N_X(v_j) &{\rm if}\,\,\,j=1\,\,\,{\rm or}\,\,\,N_X(v_j)=\emptyset.
				\end{cases}
			\end{align*}
			Then by Fact~\ref{f88}, 
			we have that 
			$N'_X(v_a)\cap N'_X(v_b)=\emptyset$ for any distinct integers $a,b\in\{1,\ldots,\ell\}.$ 
			Thus, 
			$$\left|N_X(P^*) \right|\ge \left| \bigcup_{j=1}^\ell N'_X(v_j) \right|=\sum_{j=1}^\ell |N'_X(v_j)|.$$ 
			Therefore, 
			\begin{align*}
				\sum_{i=2}^q d_{V(P^*)}(x_i)
				=|N_X(v_1)|+ \sum_{j=2}^\ell \left| N_X(v_j) \right|
				&\le |N'_X(v_1)| +\sum_{j=2}^\ell \left( |N'_X(v_j)|+1 \right)
				\le |N_X(P^*)| + \ell -1.
			\end{align*}
			This completes the proof of Fact~\ref{f7}.
		\end{proof} 
		
		\begin{fact}\label{f9}
			For each $x_i~(2\le i \le q),$ 
			$x_i$ is adjacent to at most $q$ member of $\mathcal{P}.$ 
			Hence 
			$$\sum_{P^*\in \mathcal{P}} \left| N_X(P^*) \right| \le q(q-1).$$ 
		\end{fact}
		\begin{proof}[Proof of Fact~\ref{f9}]
			Suppose to the contrary that 
			there exists a vertex $x_i$ 
			such that 
			$x_i$ is adjacent to $P^*_1,\ldots, P^*_{g},$ 
			where $g\ge q+1$ and $P^*_1,\ldots, P^*_{g}\in \mathcal{P}.$  
			For each member $P^*_b\in \mathcal{P},$ 
			choose a vertex $y_{P^*_b}\in N(x_i)\cap V(P^*_b).$ 
			We will find a $g$-fan from $x_i$ to $S\setminus\{x_i\},$  where 
			$g\ge q+1.$ 
			
			If $P^*_b=v,$ let $Q_{v}$ be the edge $x_i v.$ 
			If $P^*_b\neq v,$ write $P^*_b=\overrightarrow{P}[z,z').$ 
			Define 
			\begin{align*}
				Q_{P^*_b}
				&=\begin{cases} 
					x_i y_{P^*_b} \overleftarrow{P} z &{\rm if}\,\,\,z\in S,\\
					x_i y_{P^*_b} \overleftarrow{P} a_j \overleftarrow{P} x_j &{\rm if}\,\,\,z=a_j \in A\setminus S\,\,\,{\rm and}\,\,\,j\in\{1,\ldots, q\}\setminus\{i\},\\
					x_i y_{P^*_b} \overleftarrow{P} a_i P_i w &{\rm if}\,\,\,z=a_i \in A\setminus S.
				\end{cases}
			\end{align*} 
			The family of these internally disjoint $(x_i, S\setminus \{x_i\})$-paths $Q_{P^*_1},\ldots, Q_{P^*_g},$ is a 
			$g$-fan from $x_i$ to $S\setminus\{x_i\},$ 
			where $g\ge q+1,$ 
			which contradicts the choice of the vertex $w.$ 
			This proves Fact~\ref{f9}. 
		\end{proof}
		We consider $d_G(w)+\sum_{i=2}^{q}d_G(x_i).$
		We assert that 
		\begin{align}\label{ea1}
			\sum_{P^*\in \mathcal{P}} \left( \left| N_X(P^*) \right| -1 \right) 
			\le q(q-2).    
		\end{align}
		If fact, for each $P^*\in \mathcal{P},$ 
		we have $\left| N_X(P^*) \right|\le |X|=q-1.$ 
		If $|\mathcal{P}|\le q,$ 
		then by Fact~\ref{f2}(i), 
		\begin{align*} 
			\sum_{P^*\in \mathcal{P}} \left( \left| N_X(P^*) \right| -1 \right) \le |\mathcal{P}| (q-2)\le q(q-2),
		\end{align*}
		as desired, 
		If $|\mathcal{P}|>  q,$ then Fact~\ref{f9} implies that 
		\begin{align*} 
			\sum_{P^*\in \mathcal{P}} \left( \left| N_X(P^*) \right| -1 \right) 
			=\sum_{P^*\in \mathcal{P}}  \left| N_X(P^*) \right| -|\mathcal{P}|
			\le q(q-1)-|\mathcal{P}|
			< q(q-2). 
		\end{align*}
		Thus \eqref{ea1} holds. 
		
		By Fact~\ref{f2}(ii), we have 
		\begin{align}\label{ea3}
			\sum_{i=2}^q \sum_{j=2}^q  \left|  N_G(x_i) \cap V\left( \overrightarrow{P}[x_j,a_j ) \right)  \right|
			&=
			\sum_{i=2}^q   \left|  N_G(x_i) \cap V\left( \overrightarrow{P}[x_i,a_i ) \right)  \right|  \notag{}\\
			&\le 
			\sum_{i=2}^q \left(  \left|  \overrightarrow{P}[x_i,a_i ) \right|  -1   \right)
			= 
			\sum_{i=2}^q  \left|  \overrightarrow{P}[x_i,a_i ) \right| -q+1. 
		\end{align}
		Based on Fact~\ref{f7} and \eqref{ea1}, 
		we have 
		\begin{align*}
			\sum_{i=2}^q \sum_{P^*\in \mathcal{P}} \left| N_G(x_i) \cap V(P^*)  \right| 
			&=\sum_{P^*\in \mathcal{P}}\sum_{i=2}^q \left| N_G(x_i) \cap V(P^*)  \right| \notag{}\\
			&\le    
			\sum_{P^*\in \mathcal{P}}  
			\left( |N_X(P^*)| +|P^*| -1  \right) \notag{}\\
			&\le 
			q(q-2)
			+ \sum_{P^*\in \mathcal{P}} |P^*|.
		\end{align*} 
		This, together with \eqref{ea3}, implies 
		\begin{align*}
			\sum_{i=2}^q d_{V(P)} (x_i)
			&= 
			\sum_{i=2}^q \left| N_G (x_i) \cap V(P)\right| \\
			&= 
			\sum_{i=2}^q 
			\left( 
			\sum_{P^*\in \mathcal{P}} \left| N_G(x_i) \cap V(P^*)  \right| 
			+ 
			\sum_{j=2}^q  \left|  N_G(x_i) \cap V\left( \overrightarrow{P}[x_j,a_j ) \right)  \right|
			\right)\\
			&\le 
			q(q-2) + \sum_{P^*\in \mathcal{P}} |P^*|
			+
			\sum_{i=2}^q  \left|  \overrightarrow{P}[x_i,a_i ) \right| -q+1\\
			&=
			|P| +q^2-3q+1. 
		\end{align*}
		Combining this with Fact~\ref{f6}(i)-(ii) and Fact~\ref{f66}(iii) gives
		\begin{align*}
			d_G(w)+\sum_{i=2}^q d_G(x_i) 
			&=
			d_{V(P)} (w) +d_B (w) +d_R (w)
			+\sum_{i=2}^q \left( d_{V(P)} (x_i) +d_B(x_i) +d_R(x_i) \right) \\
			&\le 
			q+|P| +q^2-3q+1+|R|+|B|-1\\
			&=n+q^2-2q,
		\end{align*}
		which contradicts the assumption 
		that 
		$\sigma_p(S)\ge n+(p-1)^2$ for every integer $k\le  p\le \alpha(G).$ 
		This completes the proof of Theorem~\ref{thm:main}. 
	\end{proof}

	\section{Application}\label{s3cz}
	In this paper, we mainly study the degree-sum condition for the existence of a path passing through prescribed vertices with two endpoints specified.     
	Next we apply Corollary~\ref{cor31} to prove Conjecture~\ref{conj:ffgjm}. That is, 
	\begin{theorem}\label{thm:czffgjm}
		Let $G$ be a $k$-connected graph of order $n$, where $k\ge3$.  If
		\begin{align}\label{eq:cz2}
			\delta(G)>\frac{n+k^2-2k}{k}
			~~\text{and}~~
			\delta(G)\ge\alpha(G)+k-2,
		\end{align}
		then $G$ is hamiltonian-connected.
	\end{theorem} 
	
	\begin{proof}[\bf Proof of Theorem~\ref{thm:czffgjm}]
		If $\alpha(G)<k$, then the degree-sum condition in
		Corollary~\ref{cor31} is trivially satisfied, and thus
		$G$ is hamiltonian-connected, as desired.
		
		Next we assume that $\alpha(G)\ge k$, and let $k\le p\le \alpha(G).$  
		By \eqref{eq:cz2}, 
		$k\delta(G)>n+k^2-2k,$ 
		and thus 
		$$
		k\delta(G)\ge n+k^2-2k+1.
		$$
		Also, \eqref{eq:cz2} gives $\delta(G)\ge \alpha(G)+k-2\ge p+k-2$. Since $p-k\ge0$, it follows that
		\begin{align*}
			p\delta(G)
			&=(p-k)\delta(G)+k\delta(G)\\
			&\ge (p-k)(p+k-2)+n+k^2-2k+1\\
			&=n+(p-1)^2.
		\end{align*}
		Thus $\sigma_p(G)\ge p\delta(G)\ge n+(p-1)^2$ for every $k\le p\le \alpha(G).$ Then by Corollary~\ref{cor31}, $G$ is hamiltonian-connected, as desired. 
		This completes the proof of Theorem~\ref{thm:czffgjm}. 
	\end{proof} 
	
	\begin{remark}
		The condition in \eqref{eq:cz2} can not be weakened.  
		In fact, if we replace \eqref{eq:cz2} with 
		\begin{align}\label{eq:cz3}
			\delta(G)\ge\frac{n+k^2-2k}{k}
			~~\text{and}~~
			\delta(G)\ge\alpha(G)+k-2, 
		\end{align}
		then for $m\ge k-1,$ let $H:=K_k\vee kK_m.$ 
		Then $H$ satisfies \eqref{eq:cz3} and, as shown in Remark~\ref{r1}, is $k$-connected and not hamiltonian-connected.  
		On the other hand, if we replace \eqref{eq:cz2} with 
		\begin{align}\label{eq:cz4}
			\delta(G)>\frac{n+k^2-2k}{k}
			~~\text{and}~~
			\delta(G)\ge\alpha(G)+k-3, 
		\end{align} 
		then for $t\ge k+1,$ 
		let $G'=K_t\vee t K_{k-2}.$ 
		Clearly, $\alpha(G')=\kappa(G')=t\ge k+1,$ and hence $G'$ is $k$-connected.  
		Moreover, 
		$\delta(G')=t+k-3=\alpha(G') +k-3$ 
		and $k\delta(G')-(|G'|+k^2-2k)=t-k\ge 1.$ 
		Hence, $G'$ satisfies \eqref{eq:cz4}. However, $G'$ is not hamiltonian-connected, because its toughness is at most 1.  
	\end{remark}

	\section*{\normalsize Declaration}
	
	\noindent\textbf{Conflict~of~interest}
	The authors declare that they have no known competing financial interests or personal relationships that could have appeared to influence the work reported in this paper.
	
	\noindent\textbf{Data~availability}
	No data was used for the research described in the article.


\begin{thebibliography}{99}
		\bibitem{BM} J.A. Bondy, U.S.R. Murty, {\em Graph Theory,} Springer, New York, 2008.
		
		\bibitem{BLLTV}
		H.~J. Broersma, H.~Li, J.~Li, F.~Tian, and H.~J. Veldman,
		Cycles through subsets with large degree sums,
		\emph{Discrete Math.} \textbf{171} (1997), 43--54.
		
		\bibitem{CHW}
		G.~Chen, Z.~Hu, and Y.~Wu,
		A Chv\'atal--Erd\H{o}s type theorem for path-connectivity,
		\emph{Discuss. Math. Graph Theory} \textbf{44} (2024), 1247--1260.

        \bibitem{Chvatal} V. Chv\'{a}tal, Tough graphs and Hamiltonian circuits, \emph{Discrete Math.} \textbf{5}   (1973), 215--228.  

            
		\bibitem{CE}
		V.~Chv\'atal and P.~Erd\H{o}s,
		A note on Hamiltonian circuits,
		\emph{Discrete Math.} \textbf{2} (1972), 111--113.
		
		\bibitem{D}
		G.~A. Dirac,
		Some theorems on abstract graphs,
		\emph{Proc. London Math. Soc.} (3) \textbf{2} (1952), 69--81.
		
		\bibitem{GAD}
		G.~A. Dirac,
		In abstrakten Graphen vorhandene vollst\"andige $4$-Graphen und ihre
		Unterteilungen,
		\emph{Math. Nachr.} \textbf{22} (1960), 61--85.
		
		\bibitem{EGL}
		Y.~Egawa, R.~Glas, and S.~C. Locke,
		Cycles and paths through specified vertices in $k$-connected graphs,
		\emph{J. Combin. Theory Ser. B} \textbf{52} (1991), 20--29.
		
		\bibitem{FFGJM}
		J.~R. Faudree, R.~J. Faudree, R.~J. Gould, M.~S. Jacobson, and C.~Magnant,
		Chv\'atal--Erd\H{o}s type theorems,
		\emph{Discuss. Math. Graph Theory} \textbf{30} (2010), 245--256.
		
		\bibitem{FLMW}
		E.~Flandrin, H.~Li, A.~Marczyk, and M.~Wo\'zniak,
		A generalization of Dirac's theorem on cycles through $k$ vertices in
		$k$-connected graphs,
		\emph{Discrete Math.} \textbf{307} (2007), 878--884.
		
		\bibitem{G}
		R.~J. Gould,
		A look at cycles containing specified elements of a graph,
		\emph{Discrete Math.} \textbf{309} (2009), 6299--6311.
		
		\bibitem{H}
		J.~Harant,
		On paths and cycles through specified vertices,
		\emph{Discrete Math.} \textbf{286} (2004), 95--98.
		
		\bibitem{HTWEH}
		Z.~Hu, F.~Tian, B.~Wei, Y.~Egawa, and K.~Hirohata,
		Fan-type theorem for path-connectivity,
		\emph{J. Graph Theory} \textbf{39} (2002), 265--282.
		
		\bibitem{JO}
		B.~Jackson and O.~Ordaz,
		Chv\'atal--Erd\H{o}s conditions for paths and cycles in graphs and
		digraphs: a survey,
		\emph{Discrete Math.} \textbf{84} (1990), 241--254.
		
		\bibitem{M}
		K.~Menger,
		Zur allgemeinen Kurventheorie,
		\emph{Fund. Math.} \textbf{10} (1927), 96--115.
		
		\bibitem{OO}
		O.~Ore,
		Note on Hamilton circuits,
		\emph{Amer. Math. Monthly} \textbf{67} (1960), 55.
		
		\bibitem{KO}
		K.~Ota,
		Cycles through prescribed vertices with large degree sum,
		\emph{Discrete Math.} \textbf{145} (1995), 201--210.
		
		\bibitem{OY}
		K.~Ozeki and T.~Yamashita,
		A degree sum condition concerning the connectivity and the independence
		number of a graph,
		\emph{Graphs Combin.} \textbf{24} (2008), 469--483.
		
		\bibitem{W} D.B. West, Introduction to Graph Theory, Prentice Hall, Inc., 1996. 
	\end{thebibliography}
\end{document}